\documentstyle[11pt]{amsart}

\newtheorem{thm}{Theorem}
\newtheorem{prop}[thm]{Proposition}

\theoremstyle{definition}

\newtheorem{remark}[thm]{Remark}

\newcommand{\pr}{\operatorname{pr}}

\newcommand{\Proj}{\operatorname{Proj}}

\newcommand{\shom}{\operatorname{{\cal H}om}}

\newcommand{\bbP}{{\Bbb P}}

\newcommand{\cJ}{{\cal J}}
\newcommand{\cF}{{\cal F}}
\newcommand{\cO}{{\cal O}}
\newcommand{\cM}{{\cal M}}

\newcommand{\cI}{{\cal I}}

\newcommand{\cL}{{\cal L}}
\newcommand{\cT}{{\cal T}}

\newcommand{\isomo}{\overset{\sim}{=}}

 \newcommand{\brokrarr}{\vphantom{\to}\mathrel{\smash{{-}{\rightarrow}}}}

\newcommand{\PGL}{{\operatorname{PGL}}}

\newcommand{\lra}{\longrightarrow}

\let\to=\longrightarrow

\begin{document}
\title[Equivariant resolution] 
{Equivariant resolution of points of indeterminacy}
\author[Z. REICHSTEIN and B. YOUSSIN]
{Z. Reichstein and B. Youssin} 
\address{Department of Mathematics, Oregon State University 
Corvallis, OR 97331 \hfill\break
\hbox{{\rm\it\hskip\parindent Current address\/}}: 
1498 Meadows Drive, Lake Oswego, OR 97034, USA}
\thanks{Z. Reichstein was partially supported by NSF grant DMS-9801675}
\email{zinovy@@math.orst.edu}
\address{Department of Mathematics and Computer Science,
University of the Negev, Be'er Sheva', Israel\hfill\break
\hbox{{\rm\it\hskip\parindent Current address\/}}: 
Hashofar 26/3, Ma'ale Adumim, Israel}
\email{youssin@@math.bgu.ac.il}
%%%%%%%%%%%%%%%%%%%%%%%%
% 14E15 --- Birational geometry: global theory and resolution of singularities.
% 14L30 --- Group actions on varieties and schemes     
%%%%%%%%%%%%%%%%%%%%%%%

\begin{abstract}
We prove an equivariant version of Hironaka's theorem
on elimination of points of indeterminacy. Our
proof relies on canonical resolution of singularities.
\end{abstract}

\maketitle

\section{Introduction}
\label{sect1}

Throughout this note we shall work over an algebraically closed 
field of characteristic 0. All algebraic varieties, schemes, groups,
and all maps between them will be defined over $k$. The main objects 
of interest for us will be algebraic varieties 
with a $G$-action; we will refer to them as $G$-varieties. 
A $G$-equivariant morphism between two such varieties will be called
a morphism of $G$-varieties. The terms "rational map of $G$-varieties",
``birational morphism of $G$-varieties", "birational isomorphism 
of $G$-varieties", etc., are defined in a similar manner. 

Hironaka's theorem on elimination of points of indeterminacy
(see~\cite[\S0.5, Question~E and Main Theorem~II]{hironaka})
asserts that every rational 
map $f \colon X \brokrarr Y$ can be resolved into a regular map 
by a sequence of blowups $\pi \colon X_m \lra \dots \lra X_0 = X$
with smooth centers. In other words, $\pi$ can be chosen so that
the composition $f\pi$ is regular.  The purpose of this paper is to
prove the following equivariant version of this result.  

\begin{thm} \label{thm4.2}
Let $f\colon X\brokrarr Y$ be a rational map of $G$-varieties where
$Y$ is complete.
Then there is a sequence of blowups
\begin{equation} \label{tower6}
\pi\colon X_m \lra X_{m-1} \lra
\dots \lra X_1 \lra X_0 = X  
\end{equation}
with smooth $G$-invariant centers such that
the composition $f\pi$ is regular.
\end{thm}

Our proof will rely on canonical resolution of singularities. 
Along the way we prove an equivariant form 
of Chow's lemma (Proposition~\ref{cor:may28a}), generalizing 
a theorem of Sumihiro (\cite[Theorem 2]{sumihiro}).

The second author warmly thanks Institute of Mathematics of Hebrew
University for its hospitality during 1999/2000.

\section{Birational morphisms as blowups} 

The following result is an equivariant analogue 
of \cite[Theorem~7.17]{Hart}.

\begin{prop} \label{prop.b2}
Let $f\colon X'\to X$ be a birational proper morphism of $G$-varieties,
where $X$ is smooth and $X'$ is quasiprojective.
Then there exists a $G$-invariant sheaf of ideals $\cI$ on $X$ such that $X'$ 
is the blowup of $\cI$.
\end{prop}

\begin{pf}
Let $\sigma\colon G\times X\to X$ be the given action of $G$ on $X$ and
$\pr_2\colon G\times X\to X$ be the projection onto the second factor.

By a theorem of Kambayashi~\cite{kamb}, 
there is an action of $G$ on the projective space
$\bbP^n$ (via a representation $G\to \PGL_{n+1}$) 
and a $G$-equivariant embedding
$X'\hookrightarrow\bbP^n$; this yields a
$G$-equivariant embedding $i\colon X'\hookrightarrow\bbP^n\times X$.

Here $\bbP^n\times X$ is a projective space over $X$; set
$\cL=i^*\cO_{\bbP^n\times X}(1)$ and $\cT=\bigoplus_{d=0}^\infty f_*(\cL^d)$,
where $\cL^0=\cO_{\bbP^n\times X}$.
Let $\cT_1=f_*\cL$ be the component of $\cT$ of degree one.
The action of $G$ on $\bbP^n\times X$ yields a $G$-linearization of the sheaf
$\cT_1$, i.e., an isomorphism $\sigma^*\cT_1\isomo\pr_2^*\cT_1$ which satisfies
the same cocycle condition as in the definition of $G$-linearization of an
invertible sheaf (see, e.g., \cite[Definition~1.6]{git}); informally speaking,
$G$ acts on the pair $(X,\cT_1)$.

We refer to the proof of \cite[Theorem~7.17]{Hart} for the following facts:
\begin{enumerate}
\item After replacing the embedding $i$ by its $e$-fold embedding for some
positive integer $e$ (thus replacing $\cL$ by $\cL^e$), we may assume that the
graded $\cO_X$-algebra $\cT$ is generated by $\cT_1$.
\item $X'\isomo\Proj\cT$.
\item Assume $\cT$ is generated by $\cT_1$ as in (1).
If there is an invertible sheaf $\cM$ on $X$ and a sheaf of ideals $\cI$
on $X$ such that $\cI\isomo\cT_1\otimes\cM$, then $X'$ is isomorphic to the
blowup of $\cI$. 
\end{enumerate}

The variety $X$ is smooth, and hence, for any sheaf of ideals $\cF$ on $X$ of
rank one without torsion, its dual $\cF^*=\shom(\cF,\cO_X)$ is an invertible
sheaf.
(To see this, note that locally at any point $x\in X$, the generator of
$\cF^*_x$ is given by the
homomorphism $\cF_x\to\cO_{X,x}$ which maps the generators of $\cF_x$ as an
$\cO_{X,x}$-module of rank one, into elements of $\cO_{X,x}$ not having a
nontrivial common multiple; such homomorphism is unique up to an invertible
multiple since the local ring $\cO_{X,x}$ is regular, and hence, factorial.)

Thus the second dual $\cT_1^{{*}{*}}$ is an invertible sheaf, and we
have an embedding $\cT_1\hookrightarrow\cT_1^{{*}{*}}$.
The $G$-linearization of $\cT_1$ yields a $G$-linearization of
$\cT_1^{{*}{*}}$, and the above embedding is, in fact,
an embedding of $G$-linearized sheaves.
Taking $\cM=(\cT_1^{{*}{*}})^{-1}$, we see that $X'$ is isomorphic to the blowup
of the sheaf of ideals $\cI=\cT_1\otimes(\cT_1^{{*}{*}})^{-1}$.

The $G$-linearizations of $\cT_1$ and $\cT_1^{{*}{*}}$ yield a $G$-linearization
of $\cI$ and a $G$-linearized embedding
$\cI=\cT_1\otimes(\cT_1^{{*}{*}})^{-1}\hookrightarrow
\cT_1^{{*}{*}}\otimes(\cT_1^{{*}{*}})^{-1}=\cO_X$. This
shows that $\cI$ is a $G$-invariant sheaf of ideals on $X$.
\end{pf}

\section{Canonical simplification of a finite collection of ideals}

One of the main resolution theorems of
Hironaka~\cite[Main Theorem~II]{hironaka} asserts that any finite
collection of sheaves of ideals $\{\cI_i\}$ can be ``simplified''
by a finite sequence $\pi \colon X_m \lra \dots \lra X_0 = X$
of blowups with smooth centers.
In other words, the sequence
of blowups can be chosen so that $\pi^* \cI_i$ is locally principal
for each $i$, and is locally generated by a monomial with respect to a
normal crossing divisor.

Bierstone and Milman~\cite[Theorem~1.10]{bm} proved that
any sheaf of ideals $\cI$ on an algebraic variety $X$ can be
``simplified" in this sense
in a canonical way, so that the sequence $\pi$ is canonically defined;
in particular, every automorphism of $X$ 
preserving $\cI$ lifts to the entire sequence;
see~\cite[Remark~1.5]{bm}.
This immediately implies that any finite {\em ordered} collection of
sheaves of ideals $\{\cI_i\}$ can be ``simplified'' in a canonical way.

In this section we will show that a finite {\em unordered} collection 
of sheaves of ideals can be simplified in a similar manner. This result
will be used in the proof of Theorem~\ref{thm4.2}.

\begin{prop} \label{prop:may26a}
Let $X$ be an algebraic variety and $\{\cI_1,\dots,\cI_n\}$ be a finite
collection of sheaves of ideals on $X$.  
Then there exists a sequence of blowups
\begin{equation} \label{e.tower2}
 \pi\colon X_m \lra X_{m-1} \lra \dots \lra X_1 \lra X_0 = X   
\end{equation}
such that $\pi^*\cI_i$ is locally principal for
each $i$ and any automorphism of $X$ that preserves (but possibly
non-trivially permutes) the collection
$\{\cI_1, \dots, \cI_n\}$ 
lifts to the entire sequence~\eqref{e.tower2}.
\end{prop}

\begin{pf} To motivate our construction, we begin with the following 
observation.  Let $V(\cI_i)$ be the subscheme of $X$ cut out by
$\cI_i$.  If the subschemes $V(\cI_i)$ are pairwise
disjoint, i.e., if $\cI_i + \cI_j = \cO_x$ for any $i \neq j$,
then the proposition follows immediately 
from~\cite[Theorem 1.10 together with Remark~1.5]{bm}:
indeed, a sequence of blowups that simplifies 
their intersection $\cI_1 \cap \dots \cap \cI_n$, will simplify each
$\cI_i$.

The idea of the proof is to reduce the general case 
to the case where every $i$-fold intersection 
of $V(\cI_1), \dots, V(\cI_n)$ is empty (i.e. the sum of any $i$
of the sheaves $\cI_1, \dots, \cI_n$ equals $\cO_X$) first for $i = n$, then
for $i = n-1$, etc., until we reach $i = 2$.  We use descending induction
on $i$. For the base case we can take $i = n+1$, where the condition
we are interested in is trivially satisfied. 

Let $S_{\Lambda} = \sum_{j \in \Lambda} \cI_j$, where
$\Lambda$ is a subset of $\{ 1, \dots, n \}$. 
For the induction step, assume 
\begin{equation} \label{e.inters}
\text{$S_{\Lambda} = \cO_X$ whenever 
$|\Lambda| = i$,}
\end{equation}
for some $i \geq 2$.
Set $\cJ = \bigcap_{|\Omega| = i-1} S_{\Omega}$.
Note that by our assumption $S_{\Omega_1} + S_{\Omega_2} = \cO_{X}$
for any two distinct subsets $\Omega_1$ and $\Omega_2$ of
$\{1, \dots, n \}$ of cardinality $i-1$, so that $V(\cJ)$ is the
disjoint union of $V(S_{\Omega})=\bigcap_{j\in\Omega}V(\cI_j)$ with
$|\Omega| = i-1$.

Let $\pi \colon X' \lra X$ be the canonical 
simplification of the sheaf $\cJ$; as we have seen above, $\pi$
simplifies each $S_{\Omega}$ with $|\Omega| = i-1$.
For $j = 1, \dots, n$, denote the conductor
$(\pi^* \cI_j):(\pi^* \cJ)$ by $\cI'_j$; it is natural to think of
$\cI'_j$ as a ``weak transform" of $\cI_j$.
The stalk of this sheaf of ideals at a (not necessarily closed) point
$x \in X'$ is described as follows.  

If $x\notin V(\pi^* \cI_j)$ then
$(\cI_j')_x=(\pi^* \cI_j)_x=\cO_{x,X'}$.
If $x\in V(\pi^* \cI_j)$ and $x\in V(\pi^*S_{\Omega})$ for some
$\Omega$
satisfying $|\Omega| = i-1$, then such $\Omega$ is unique
and $j \in \Omega$ (otherwise 
$V(\pi^*S_{\Omega\cup\{j\}})$ would be nonempty, contrary 
to~\eqref{e.inters}). 
Thus in this case $\pi^* \cI_j\subset \pi^*S_{\Omega}$ and
$(\cI_j')_x=(\pi^* S_{\Omega})_x^{-1}(\pi^* \cI_j)_x$, where
$(\pi^* S_{\Omega})_x\subset\cO_{x,X'}$ is a principal ideal.
Finally, if $x\in V(\pi^* \cI_j)$ and $x\notin V(\pi^*S_{\Omega})$ for any
$\Omega$ satisfying $|\Omega| = i-1$, then $x\notin V(\pi^*\cJ)$ and
$(\cI_j')_x=(\pi^* \cI_j)_x$.
To summarize:
\begin{equation} \label{eqn:jun6a}
(\cI_j')_x=
\begin{cases}
\ (\pi^* S_{\Omega})_x^{-1}(\pi^* \cI_j)_x&
\text{ if $x\in V(\pi^*S_{\Omega})$ for some $\Omega\subset 
\{ 1, \dots , n \} $} \\ 
 &      \text{ such that $j \in \Omega$ and
$|\Omega| = i-1$,}\\
\ (\pi^* \cI_j)_x &
\text{ otherwise.}
\end{cases}
\end{equation}
Consequently, for any sequence of blowups $\pi'\colon X''\lra X'$,
the ideal $(\pi')^*\cI_j'$ is locally principal if and only if the
ideal $(\pi')^*(\pi^*\cI_j)$ is locally principal. 
This reduces the problem of simplifying the collection
$\{ \cI_1, \cdots, \cI_n \}$ of sheaves of ideals on $X$ to
the problem of
simplifying the collection $\{ \cI_1', \cdots, \cI_n' \}$ of sheaves
of ideals on $X'$.

For $\Lambda\subset\{1,\dots,n\}$, set
$S'_{\Lambda} = \sum_{j \in \Lambda} \cI'_j$.
We claim that
\begin{equation} \label{e.inters'}
 \text{$S'_{\Lambda} = \cO_{X'}$ whenever 
$|\Lambda| = i-1$.} 
\end{equation}
We will prove this equality by showing that 
$(S'_{\Lambda})_x = \cO_{X',x}$ for every $x \in X'$. 
Indeed, if $x \not \in V(\pi^*(\cI_j))$ for some $j \in \Lambda$
then \[ \cO_{X', x} = (\pi^* \cI_j)_x \subset (\cI'_j)_x \subset 
(S'_{\Lambda})_x \, , \]
as desired. On the other hand, if $x \in V(\pi^* \cI_j)$ for 
every $j \in \Lambda$, i.e., $x \in V(\pi^* S_{\Lambda})$, 
then~\eqref{eqn:jun6a} tells us that
\[ (S'_{\Lambda})_x = \sum_{j \in \Lambda} (\cI'_j)_x =
\sum_{j \in \Lambda} (\pi^* S_{\Lambda})_x^{-1}(\pi^*\cI_j)_x  = 
(\pi^* S_{\Lambda})_x^{-1} (\pi^* S_{\Lambda})_x = \cO_{X', x} \, . \]  

We have thus reduced the problem
of simplifying the collection $\{ \cI_1, \cdots, \cI_n \}$ of sheaves
of ideals on $X$, satisfying condition~\eqref{e.inters},
to the problem of simplifying the collection
$\{ \cI'_1, \dots, \cI'_n \}$ of sheaves of ideals on $X'$, satisfying
condition~\eqref{e.inters'}. This completes the induction step.

To finish the proof of the proposition, note that
the sequence~\eqref{e.tower2} of blowups 
constructed by the recursive algorithm we just described, 
depends on $X$ and the unordered collection
$\{\cI_i\}$ in a canonical way; see~\cite[Remark~1.5]{bm}. Hence,
any automorphism of $X$ that preserves the unordered collection
$\{\cI_i\}$, lifts to to the entire sequence~\eqref{e.tower2}, as claimed.
\end{pf}

\begin{remark} \label{rem.simpl} Our proof also shows that 
each $\pi^*{\cI}_i$ is generated by a monomial with respect to
a normal crossing divisor. In other words, $\pi$ simplifies each
$\cI_i$ in the sense of Hironaka's original definition; for details
see~\cite[Remark~1.8]{bm}.
This assertion will not be used in the sequel; for this reason
we did not include it in the statement of Proposition~\ref{prop:may26a}.
\end{remark}

\section{Equivariant Chow Lemma}  

In this section we will prove the following generalization of
Chow's lemma. 

\begin{prop}
\label{cor:may28a}
For every $G$-variety $X$, there exists a quasiprojective $G$-variety
$Z$ and a proper birational morphism $Z \to X$. If $X$ is complete
then $Z$ is projective. 
\end{prop}

Note that if $G$ is assumed to be connected, this result 
is a well-known theorem of Sumihiro~\cite[Theorem~2]{sumihiro}; 
see also~\cite[Theorem~1.3]{pv}. The argument below reduces the general
case to the case where $G$ is connected.

\begin{pf} The second assertion is an immediate consequence of the first:
if $X$ is complete, $Z \lra X$ is proper and $Z$ is quasiprojective 
then $Z$ is also complete and, hence, projective. 

To prove the first assertion
let $G_0$ be the connected component of $G$.
Applying Sumihiro's theorem to $X$, viewed as a $G_0$-variety,
yields a quasiprojective $G_0$-variety $Y$ and a proper
$G_0$-equivariant birational morphism $f\colon Y\to X$.

For the rest of this proof we shall use set-theoretic notation:
by a ``point" we will always mean a closed point. 

Recall that the homogeneous fiber product $G *_{G_0} Y$ is the $G$-variety
defined as the geometric
quotient $(G\times Y)/G_0$ for the action of $G_0$ given by
$g_0(g,y)=(gg_0^{-1},g_0y)$, where $g\in G$, $g_0\in G_0$ and $y\in Y$;
see~\cite[Section 4.8]{pv}.  We shall write $[g, y]$ for the element of
$G*_{G_0} Y$ represented by $(g, y) \in G \times Y$. Since $G_0$
has finite index in $G$, $G *_{G_0} Y$ admits a more concrete description
as a disjoint union of $|G/G_0|$ copies of $Y$. More precisely,
if we choose a representative $a_h$ for each $h \in G/G_0$, we can explicitly
identify $G/G_0 \times Y$ and $G *_{G_0} Y$ as abstract varieties,
via $(h, y) \mapsto [a_h, y]$. Moreover, if we define a $G$-action
on $G/G_0 \times Y$ by $g (h, y) \lra (\overline{g}h, 
(a_h^{-1} g^{-1} a_{\overline{g}h}) y)$ then $(h, y) \mapsto [a_h, y]$ 
identifies $G/G_0 \times Y$ and $G *_{G_0} Y$ as $G$-varieties.
Here $\overline{g}$ is the image of $g$ in $G/G_0$
and $(a_h^{-1} g^{-1} a_{\overline{g}h})y$ is well-defined because
$a_h^{-1} g^{-1} a_{\overline{g}h}$ is an element of $G_0$. 

Let $\alpha \colon G *_{G_0} Y \lra X$ and $\beta \colon G *_{G_0} Y 
\lra G/G_0$ be the maps of $G$-varieties given by 
$\alpha \colon [g, y] \lra gf(y)$ and
$\beta \colon [g, y] \mapsto \overline{g}$. (Here $G$ acts on $G/G_0$
by left multiplication.) These maps are shown in the diagram below.
\[ \begin{array}{ccccc}
      &                  & G *_{G_0} Y \simeq G/G_0 \times Y &      & \\
      & \alpha \swarrow  &   &  \beta \searrow \nwarrow s &  \\
  X    &                  &                         &           &  G/G_0  
\end{array} \]
Let $S$ be the set of all sections $s$ of $\beta$. Note that
if we identify $G *_{G_0} Y$ with $G/G_0 \times Y$ as above, then
$\beta \colon G /G_0 \times Y \lra G/G_0$ is the projection 
to the first factor. Thus $S \simeq Y^{|G/G_0|}$  
as an abstract variety.  Moreover, since $\beta$ is $G$-equivariant,
$G$ acts on this variety by $g \colon s \mapsto t$,
where $s,t \in S$ and $t(h) = g \cdot s(\overline{g}^{-1}h)$ 
for any $h \in G/G_0$. 

Let $Z$ be the closed $G$-invariant subvariety of
$S$ consisting of those sections $s \colon G/G_0 \lra G *_{G_0} Y 
\simeq G/G_0 \times Y$ with the property that $\alpha \circ s(G/G_0)$ 
is a single point of $X$; we shall denote this point by $x_s$. 

We claim that the morphism $\phi \colon Z \lra X$ given by $s \lra x_s$
has the properties asserted in the proposition.
Indeed, since $Y$ is quasiprojective, 
and $Z$ is a closed subvariety of $S \simeq Y^{|G/G_0|}$, 
$Z$ is quasiprojective as well.  

To show that $\phi$ is a birational morphism, assume
the birational morphism $f \colon Y \lra X$ is an isomorphism 
over a dense open subset $U \subset X$. Then 
$V = \bigcap_{h \in G/G_0} a_h U$ is also a dense open subset of $X$, 
and for every $x \in V$, 
$ \alpha^{-1}(x) = \{ [a_h, f^{-1}(a_{h^{-1}}x)] \, : \, 
h \in G/G_0 \}$;
it is the image of the unique section $s_x\in S$ satisfying
$x_{s_x}=x$.
This section is given by
$s_x(h) = [a_h, f^{-1}(a_h^{-1} x)]$,
and the morphism $V \lra Z$, $x\mapsto s_x$, is a two-sided rational
inverse to $\phi$. 
\end{pf}

\section{Proof of Theorem~\ref{thm4.2}}
We begin with two reductions. First of all,
we may assume without loss of generality that
$Y$ is projective. Indeed, Proposition~\ref{cor:may28a} 
yields a projective $G$-variety $Y'$ and a
$G$-equivariant birational morphism $u\colon Y'\to Y$.
We can now replace $Y$ by $Y'$ and $f$ by
$f' = u^{-1}f\colon X\brokrarr Y'$
If we can construct a sequence of blowups $\pi$, as in~\eqref{tower6},
so that $f'\pi$ is regular, then
$f\pi$ is regular as well, i.e., the same sequence 
of blow ups will resolve the indeterminacy locus of $f$.

Secondly, we may assume that $X$ is smooth.
Indeed, let
\begin{equation} \label{eqn14mar2}
X_l@>\pi_l>>\dots @>\pi_1>>X_0=X\ ,
\end{equation}
be the canonical resolution of singularities of $X$,
as in~\cite[Theorem~7.6.1]{Vil2} or~\cite[Theorem~13.2]{bm}. Here 
$X_l$ is smooth, the centers $C_i\subset X_{i}$ are smooth 
and $G$-invariant, and the action of $G$ lifts to
the entire resolution sequence \eqref{eqn14mar2}.
Replacing $X$ by $X_l$, we may assume that $X$ is smooth.

 From now on we will assume $X$ is smooth and $Y$ is projective.
Let $G_0$ be the connected component of $G$.
As $X$ is smooth, it is normal, and we can apply
a theorem of Sumihiro (see~\cite[Theorem~1]{sumihiro} or
~\cite[Theorem 1.1]{knop-et-al} or \cite[Theorem~1.2]{pv}) 
which yields a finite covering
$X=U_1 \cup \dots \cup U_d$, where each $U_i$ is a $G_0$-invariant 
quasiprojective open subvariety of $X$.  

For each $i$, let $Z_i$ be the closure of the graph of
$f|_{U_i}\colon U_i\brokrarr Y$ in $U_i\times Y$; it is a
quasiprojective variety, and the projection $h_i\colon Z_i\to U_i$ is
a proper birational morphism of $G$-varieties.
By Proposition~\ref{prop.b2}, we can find a $G_0$-invariant
sheaf of ideals $\cI_i$ on $U_i$
such that $Z_i$ is isomorphic to the blowup of $\cI_i$.
Let $\cI_i'$ be the maximal sheaf of ideals on $X$ such that
$\cI_i'|_{U_i}=\cI_i$; as $\cI_i'$ is unique, it is $G_0$-invariant.
The $G$-invariant collection of sheaves $\{g^*\cI_i' \, |\, g \in G, \, 
i = 1, \dots, d \}$ on $X$ is finite (it contains
no more than $|G/G_0|$ sheaves for each $i$).
Proposition~\ref{prop:may26a} yields a $G$-equivariant 
sequence of blowups
\[ \pi\colon X_m \lra X_{m-1} \lra \dots \lra X_1 \lra X_0 = X \]  
with the property that the pullback $\pi^*\cI_i'$ is locally principal
for each $i$; hence, the composition
$h_i^{-1}\pi\colon X_m\to X\brokrarr Z_i$ is regular on $\pi^{-1}(U_i)$.
This implies that $f\pi\colon X_m\to X\brokrarr Z_i\to Y$ is regular on
$\pi^{-1}(U_i)$ for each $i$, and hence, on all of $X_m$.
\qed

\begin{remark} \label{rem.final} Note that 
if $Y$ is assumed to be projective in the statement of 
Theorem~\ref{thm4.2}
then Proposition~\ref{cor:may28a} is not needed in the proof. 
On the other hand, if $G$ is assumed to be connected,
then Proposition~\ref{prop:may26a} may be replaced 
by~\cite[Theorem~1.10 together with Remark~1.5]{bm}
and~Proposition~\ref{cor:may28a} may be replaced 
by~\cite[Theorem~2]{sumihiro}.
\end{remark}


\begin{thebibliography}{ABCD}

\bibitem[BM]{bm}E. Bierstone, P. D. Milman, {\em Canonical 
desingularization in characteristic zero by blowing up the maximum 
strata of a local invariant}, Invent. math. {\bf 128} (1997), no. 2, 207--302.

\bibitem[Ha]{Hart} R. Hartshorne.  Algebraic geometry.  Springer, 1977.

\bibitem[Hi]{hironaka}H. Hironaka, {\em Resolution of singularities 
of an algebraic variety over a field of characteristic 0, I and II},
Annals of Math. {\bf 79} (1964), 109--326.

\bibitem[Ka]{kamb}
T. Kambayashi, {\em Projective representations of algebraic groups of
transformations}, Amer. J. Math. {\bf 88} (1966), 199--205.

\bibitem[KKLV]{knop-et-al}F. Knop, H. Kraft, D. Luna, Th. Vust,
{\em Local properties of algebraic group actions}, in ``Algebraische
Transformationsgruppen und Invariantentheorie", DMV {\bf 13}, Birkh\"auser,
Basel, 1989.

\bibitem[MFK]{git}D. Mumford, J. Fogarty and F. Kirwan.
Geometric invariant theory.  Third enlarged edition, Springer, 1994.

\bibitem[PV]{pv}V. L. Popov, E. B. Vinberg, {\em Invariant Theory}, in
Encyclopedia of Math. Sciences {\bf 55}, Algebraic Geometry IV, edited by A. N. Parshin and I. R. Shafarevich, Springer-Verlag, 1994. 

\bibitem[Su]{sumihiro} H. Sumihiro, {\em Equivariant completion},
J. Math. Kyoto Univ., {\bf 14}, no. 1 (1974), pp.~1--28

\bibitem[V]{Vil2} O. E. Villamayor U., {\em Patching local uniformizations,}
Ann. scient. \'Ec. Norm. Sup., 4e s\'erie, {\bf 25} (1992), 629--677.
\end{thebibliography}
\end{document}